\documentclass[11pt]{article}
\usepackage[table]{xcolor}
\usepackage{graphicx}
\usepackage{amsmath,amsfonts,amssymb}
\usepackage{graphicx}
\usepackage{multirow}
\usepackage{tikz,pgf}
\usepackage{url}
\usepackage{amsmath}
\usepackage{graphicx}
\usepackage{epsfig}
\usepackage{epstopdf}
\usepackage{color}
\usepackage{enumitem}
\def\disp{\displaystyle}

\def\tto{\;{\lower 1pt \hbox{$\rightarrow$}}\kern -10pt
\hbox{\raise 2pt \hbox{$\rightarrow$}}\;}

\def\Bar{\overline}

\def\epsilon{\varepsilon}

\def\h{\hfill\Box}
\def\R{\Bbb R}

\def\ox{\bar{x}}

\def\oy{\bar{y}}
\def\oz{\bar{z}}

\def\co{\mbox{\rm core}}
\def\core{\mbox{\rm core}\,}

\def\gph{\mbox{\rm gph}\,}

\def\epi{\mbox{\rm epi}\,}

\def\dom{\mbox{\rm dom}\,}

\def\ker{\mbox{\rm ker}\,}

\def\h{\hfill\square}
\def\dn{\downarrow}
\def\ph{\varphi}
\def\emp{\emptyset}

\def\oR{\Bar{\R}}

\def\gg{\gamma}
\def\dd{\delta}

\def\ph{\varphi}
\def\emp{\emptyset}

\def\oR{\Bar{\R}}

\def\gg{\gamma}
\def\dd{\delta}

\setlist[enumerate,1]{itemsep=0.0ex,parsep=0.5ex,label={\rm(\alph*)},leftmargin=*,align=left}
\newcounter{lk}

\usepackage[colorlinks=true]{hyperref}

\begin{document}
\begin{center}
{\sc\bf ALGEBRAIC CORE AND CONVEX CALCULUS WITHOUT TOPOLOGY}\\[1ex]
Dang Van Cuong\footnote{Department of Mathematics, Faculty of Natural Sciences, Duy Tan University, Da Nang, Vietnam (dvcuong@duytan.edu.vn)}, Boris S. Mordukhovich\footnote{Department of Mathematics, Wayne State University, Detroit, Michigan 48202, USA (boris@math.wayne.edu). Research of this author was partly supported by the USA National Science Foundation under grants DMS-1512846 and DMS-1808978, by the USA Air Force Office of Scientific Research grant \#15RT04, and by the Australian Research Council under Discovery Project DP-190100555.}, Nguyen Mau Nam\footnote{Fariborz Maseeh Department of Mathematics and Statistics, Portland State University, Portland, OR 97207, USA (mnn3@pdx.edu). Research of this author was partly supported by the USA National Science Foundation under grant DMS-1716057.}, Addison Cartmell\footnote{Fariborz Maseeh Department of
 Mathematics and Statistics, Portland State University, Portland, OR 97207, USA (addison3@pdx.edu).}\\[2ex]
{\bf Dedicated to Alfredo Iusem on the occasion of his 70th birthday}
\end{center}
\small{\bf Abstract.} In this paper we study the concept of algebraic core for convex sets in general vector spaces without any topological structure and then present its applications to problems of convex analysis and optimization. Deriving  the equivalence between the Hahn-Banach theorem and and a simple version of the separation theorem of convex sets in vector spaces allows us to develop a geometric approach to generalized differential calculus for convex sets, set-valued mappings, and extended-real-valued functions with qualification conditions formulated in terms of algebraic cores for such objects. We also obtain a precise formula for computing the subdifferential of optimal value functions associated with convex problems of parametric optimization in vector spaces. Functions of this type play a crucial role in many aspects of convex optimization and its applications.\\[1ex]

\newtheorem{Theorem}{Theorem}[section]
\newtheorem{Proposition}[Theorem]{Proposition}
\newtheorem{Remark}[Theorem]{Remark}
\newtheorem{Lemma}[Theorem]{Lemma}
\newtheorem{Corollary}[Theorem]{Corollary}
\newtheorem{Definition}[Theorem]{Definition}
\newtheorem{Example}[Theorem]{Example}
\renewcommand{\theequation}{\thesection.\arabic{equation}}
\normalsize

\section{Introduction}
\setcounter{equation}{0}

Convex analysis and its numerous applications in infinite-dimensional spaces have been largely developed under certain interiority assumptions on convex sets and related objects in topological spaces; see, e.g., the books \cite{Bauschke2011,Borwein2000,bot,bi,hol,ktz,pr,z} and the references therein. Various notions of convergence of sets and functions play a prominent role in developing important results of convex analysis and applications.

Since conventional interiority conditions (involving nonempty interiors of convex sets) fail to fulfill for important classes of infinite-dimensional problems in optimization and economic modeling, more relaxed notions have been studied and applied in infinite dimensions. Concerning convex sets, these notions include various relative interior and core constructions; see \cite{bao-mor,BG,bl,bot,fl,hs,hol,ktz,z} among other publications. Most of them employ advantages of a topological structure on the space in question, while pure algebraic constructions have been also partly investigated and applied.

In this paper we concentrate on the concept of {\em algebraic core} for convex sets defined in an arbitrary (real) vector space without any topology. Algebraic cores and related nontopological notions have modestly used in the literature on nonlinear analysis and optimization, particularly in applications to vector and set-valued optimization; see, e.g., \cite{hjn,hol,ktz,z}. However, broader applications require developing {\em generalized differential calculus} for convex sets, set-valued mappings, and extended-real-valued functions in vector spaces under qualification conditions expressed in terms of algebraic cores. The main goal of this paper is to develop such a calculus, together with other useful results involving algebraic cores.

Note that it has been realized in convex analysis (starting with the finite-dimensional framework of \cite{r}) that the closedness of the sets in question is not needed for deriving basic calculus rules, although is it required for other important results. This is a striking difference from general variational analysis, where the closedness of sets and lower semicontinuity of functions is needed everywhere; see, e.g., \cite{Borwein2000,m-book1,m-book,rw}. It is due to fact that variational techniques are based on perturbation and approximation procedures with the subsequent passage to the limit, while convex analysis does not require this on a regular basis. However, a certain topological structure is an essential framework for known results and proofs in generalized differential calculus of convex analysis in infinite dimensions; see, e.g., the books \cite{Bauschke2011,Borwein2000,bot,pr,z} among many other publications.

To reach our goal on developing convex generalized differential calculus without topology, we employ a {\em dual-space geometric approach} to deal first with normals to set intersections and then apply it to deriving basic calculus rules for coderivatives of set-valued mappings and subgradients of nonsmooth functions. This approach is borrowed from variational analysis \cite{m-book1,m-book}, where it is based on the {\em extremal principle} for systems of closed sets. An appropriate version of the extremal principle for convex sets \cite{mor-nam} does not require closedness, but the topological structure is essential. Furthermore, it is shown in \cite{mor-nam} and \cite{mnrt} that the convex extremal principle is equivalent to {\em convex separation} of sets under certain {\em interiority} conditions in normed and linear convex topological vector spaces, respectively.

In this paper we rely on a proper version of convex separation theorem, which holds in any vector space and is formulated via {\em algebraic core conditions} instead of the conventional interiority assumptions in topological settings. We show also that a simple ``extreme" version of this result is equivalent the (analytic) Hahn-Banach theorem in vector spaces. Furthermore, to proceed with deriving major convex calculus rules, we need a vector space counterpart of Rockafellar's finite-dimensional result on relative interiors of convex graphs, which is obtained here in terms of algebraic cores.

The rest of the paper is organized as follows. After presenting basic definitions, we collect in Section~\ref{prelim} those properties of algebraic cores that are used below. Section~\ref{separ} revolves around separation of convex sets without topology. We provide several versions of convex separation in terms of the algebraic core and show that one of them, which constitutes an {\em extreme case of separation}, implies the Hahn-Banach extension theorem in vector spaces. This section also presents an algebraic vector space counterpart, in terms of algebraic cores and algebraic closures, of a fundamental result of finite-dimensional geometry involving relative interiors and topological closures of convex sets. The subsequent Section~\ref{graph} establishes, with the usage of convex separation, a precise {\em core representation} for {\em graphs} of convex set-valued mappings between vector spaces.

In Section~\ref{normal} we use the separation technique and core properties to derive the basic {\em intersection formula} for {\em normals} to convex sets in vector spaces under the core qualification condition. This result is employed in Section~\ref{cod-sum} and Section~\ref{cod-chain} to obtain {\em sum} and {\em chain rules}, respectively, for {\em coderivatives} of convex-graph set-valued mappings between vector spaces. The obtained results for coderivatives yield the corresponding calculus rules for {\em subgradients} of extended-real-valued convex functions under appropriate qualification conditions in terms of algebraic cores. Finally, in Section~\ref{marginal} we use algebraic cores to give a precise calculation of subgradients for {\em optimal value/marginal functions} in vector spaces, which play a crucial role in many aspects of constrained optimization and applications.

Note that a similar approach would allow us to derive generalized differential calculus rules for convex objects in locally convex topological vector spaces by using interior qualification conditions and their modifications instead of those established in this paper via algebraic cores. On the other hand, it is possible to develop a converse approach by using the strongest locally convex topology on the vector spaces in question. We prefer here a direct core algebraic approach, which does not rely on any topology.

Throughout this paper we employ the conventional notation of convex and variational analysis; see, e.g., \cite{m-book1,rw,z}. All the spaces under consideration are real vector spaces. Given such a space $X$, its {\em algebraic dual} space is defined by
\begin{equation*}
X^\prime:=\big\{f\colon X\to\R\;\big|\;f\;\text{ is a linear function}\big\}.
\end{equation*}

\section{Basic Definitions and Some Algebraic Properties}\label{prelim}
\setcounter{equation}{0}

Let us start with the basic constructions used in this paper for arbitrary vector spaces $X$. Given a nonempty set $\Omega\subset X$, define the {\em algebraic core} of $\Omega$ by
\begin{eqnarray}\label{core}
\mbox{\rm core}(\Omega):=\big\{x\in\Omega\;\big|\;\forall v\in X,\;\exists\delta>0,\;\forall t\mbox{ with }|t|<\delta:\;x+t v\in\Omega\big\}.
\end{eqnarray}
Algebraic cores are also known in the literature as ``algebraic interiors" of convex sets. A complementary notion is called the {\em algebraic closure} of $\Omega$ and is defined by
\begin{eqnarray}\label{lin}
\mbox{\rm lin}(\Omega):=\big\{x\in X\;\big|\;\exists w\in\Omega:\;[w,x)\subset\Omega\big\}.
\end{eqnarray}
Note that $[w,w)=\{w\}$. When $X$ is a topological vector space, it is easy to check the validity of the following inclusions, which all may be strict:
\begin{eqnarray*}
\mbox{\rm int}(\Omega)\subset\mbox{\rm core}(\Omega)\subset\Omega\subset\mbox{\rm lin}(\Omega)\subset\Bar{\Omega},
\end{eqnarray*}
where ${\rm int}(\Omega)$ and $\Bar{\Omega}$ signify the (topological) interior and closure of $\Omega$, respectively. Recall that a subset $\Omega$ of a vector space $X$ is {\em absorbing} if for any $v\in X$ there exists $\delta>0$ such that $t v\in\Omega$ whenever $|t|<\delta$. It follows directly from these definitions that $\ox\in\mbox{\rm core}(\Omega)$ if and only if the shifted set  $\Omega-\ox$ is absorbing. Observe also the following useful formula for representing algebraic cores of {\em set products} in vector spaces:
\begin{equation}\label{prod}
\co(\Omega\times\Omega_2)=\co(\Omega_1)\times\co(\Omega_2).
\end{equation}

For the reader's convenience, we further collect in this section some elementary properties of algebraic cores and closures of sets in vector spaces that are used in what follows.\vspace*{-0.1in}

\begin{Proposition}\label{core-conv} Let $\Omega$ be a convex subset of $X$. Then the sets $\mbox{\rm core}(\Omega)$ and $\mbox{\rm lin}(\Omega)$ are also convex in this space.
\end{Proposition}\vspace*{-0.1in}
{\bf Proof.} Fix any $a,b\in\mbox{\rm core}(\Omega)$ and $0<\lambda<1$. It follows from definition \eqref{core} that for any $v\in X$ there exists $\delta>0$ such that
\begin{equation*}
a+t v\in\Omega\;\mbox{ and }\;b+t v\in\Omega\;\mbox{ whenever }\;|t|<\delta.
\end{equation*}
Using the convexity of $\Omega$, for each such number $t$ we have
\begin{equation*}
\lambda a+(1-\lambda)b+t v=\lambda(a+t v)+(1-\lambda)(b+t v)\in\lambda\Omega+(1-\lambda)\Omega\subset\Omega.
\end{equation*}
It implies that $\lambda a+(1-\lambda)b\in\mbox{\rm core}(\Omega)$, and hence $\mbox{\rm core}(\Omega)$ is convex.

To proceed with the verification of convexity for the algebraic closure, pick any vectors $a,b\in\mbox{\rm lin}(\Omega)$ and $0<\lambda<1$. Then there exist vectors $u,v\in\Omega$ such that
\begin{equation*}
[u,a)\subset\Omega\;\mbox{ and }\;[v,b)\subset\Omega.
\end{equation*}
Denoting $x_\lambda:=\lambda a+(1-\lambda)b$ and $w_\lambda:=\lambda u+(1-\lambda)v\in\Omega$, we see that $[w_\lambda,x_\lambda)\subset\Omega$, and so $x_\lambda\in\mbox{\rm lin}(\Omega)$. This verifies the convexity of $\mbox{\rm lin}(\Omega)$. $\h$\vspace*{-0.1in}

\begin{Proposition}\label{aip} Let $\Omega\subset X$ be convex. If $a\in\mbox{\rm core}(\Omega)$ and $b\in\Omega$, then $[a,b)\subset\mbox{\rm core}(\Omega)$.
\end{Proposition}\vspace*{-0.1in}
{\bf Proof.} Fix $\lambda\in(0,1)$, define $x_\lambda:=\lambda a+(1-\lambda)b$, and then verify that $x_\lambda\in\mbox{\rm core}(\Omega)$. Since $a\in\mbox{\rm core}(\Omega)$, for any $v\in X$ there exists $\delta>0$ such that
\begin{equation*}
a+tv\in\Omega\;\mbox{ whenever }\;|t|<\delta.
\end{equation*}
Now taking such $t$ and using the convexity of $\Omega$ readily imply that
\begin{equation*}
x_\lambda+ t\lambda v=\lambda a+(1-\lambda)b+t\lambda v=\lambda(a+tv)+(1-\lambda)b\in\Omega,
\end{equation*}
which amount to saying that $x_\lambda\in\mbox{\rm core}(\Omega)$. $\h$\vspace*{-0.1in}

\begin{Proposition}\label{Lm-co-positive} Let $\Omega\subset X$ be convex, and let $x_0\in\Omega$. Suppose further that for any $v\in X$ there exists a number $\delta>0$ such that whenever $0<\lambda<\delta$ we have $x_0+\lambda v\in\Omega$. Then $x_0\in\mbox{\rm core}(\Omega)$.
\end{Proposition}\vspace*{-0.1in}
{\bf Proof.} Fix any $v\in X$ and find $\delta_+>0$ satisfying
\begin{equation*}
x_0+\lambda v\in\Omega\;\mbox{\rm for all }\;0<\lambda <\delta_+.
\end{equation*}
Furthermore, there exists a positive number $\delta_-$ such that $x_0+\lambda(-v)\in\Omega$ whenever $0<\lambda<\delta_{-}$. Letting $\delta:=\min\{\delta_+,\delta_-\}>0$, we can easily see that $x_0+\lambda v\in\Omega$ whenever $|\lambda|<\delta$. It shows that $x_0\in\mbox{\rm core}(\Omega)$. $\h$\vspace*{-0.1in}

\begin{Proposition}\label{2core} Let $\Omega$ be a convex subset of $X$. Then we have
\begin{equation*}
\mbox{\rm core}\big(\mbox{\rm core}(\Omega)\big)=\mbox{\rm core}(\Omega).
\end{equation*}
\end{Proposition}\vspace*{-0.1in}
{\bf Proof.} Note first that the set $\mbox{\rm core}(\Omega)$ is convex by Proposition~\ref{core-conv}, and hence the set ${\rm core}(\mbox{\rm core}(\Omega))$ is also convex. Since $\mbox{\rm core}(\Omega)\subset\Omega$, it follows that
\begin{equation*}
\mbox{\rm core}\big(\mbox{\rm core}(\Omega)\big)\subset\mbox{\rm core}(\Omega).
\end{equation*}
To verify the opposite inclusion, fix $a\in\mbox{\rm core}(\Omega)$ and take any $v\in X$. It follows from the definition that there exists $\delta>0$ such that
\begin{equation*}
a+tv\in\Omega\;\mbox{ whenever }\;|t|<\delta,
\end{equation*}
and hence $a+\frac{\delta}{2}v\in\Omega$. For any $\gamma$ with $0<\gamma<\frac{\delta}{2}$ define the number
\begin{equation*}
\lambda:=1-\frac{2\gamma}{\delta}.
\end{equation*}
Since $0<\gamma<\frac{\delta}{2}$, we get $\lambda\in(0,1)$, and thus Proposition~\ref{aip} tells us that
\begin{equation*}
a+\gamma v=\lambda a+(1-\lambda)\Big(a+\frac{\delta}{2}v\Big)\in\core(\Omega)\;\text{ for all }\;\gamma\;\text{ with }\;0<\gamma<\frac{\delta}{2}.
\end{equation*}
Employing now Proposition~\ref{Lm-co-positive} yields $a\in\mbox{\rm core}(\mbox{\rm core}(\Omega))$. $\h$\vspace*{-0.1in}

\begin{Proposition}\label{lm-Core-equa} Let $\Omega$ be a convex subset of $X$ with $\mbox{\rm core}(\Omega)=\Omega$. Then for any set $A \subset X$, we have the equality
$$
\mbox{\rm core}(\Omega+A)=\Omega+A.
$$
\end{Proposition}\vspace*{-0.1in}
{\bf Proof.} Observe that
\begin{align*}
\Omega+A&=\bigcup_{a\in A}(\Omega+a)=\bigcup_{a\in A}\big(\text{core}(\Omega)+a\big)=\bigcup_{a\in A}\big(\text{core}(\Omega+a)\big)\subset\text{core}(\Omega+ A).
\end{align*}
Since the opposite inclusion is obvious, the conclusion of the proposition follows. $\h$\vspace*{-0.1in}

\begin{Proposition}\label{lm1} Let $\Omega$ be a subset of $X$ with $\mbox{\rm core}(\Omega)\ne\emptyset$, and let $f\colon X\to\R$ be a nonzero linear function. Then $f$ cannot be a constant function on $\Omega$.
\end{Proposition}\vspace*{-0.1in}
{\bf Proof.} Arguing by contradiction, suppose that
\begin{equation*}
f(x)=c\;\mbox{\rm for all }\;x\in\Omega
\end{equation*}
for some constant $c$. Fix $x_0\in\mbox{\rm core}(\Omega)$ and let $\Theta:=\Omega-x_0$. Then $0\in\mbox{\rm core}(\Theta)$ and
\begin{equation*}
f(x)=0\;\mbox{\rm for all }\;x\in\Theta,
\end{equation*}
Taking any $v\in X$ and choosing $t>0$ sufficiently small such that $tv\in\Theta$ give us $f(tv)=tf(v)=0$, and thus we get $f(v)=0$. $\h$\vspace*{-0.1in}

\section{Convex Separation and Consequences in Vector Spaces}\label{separ}
\setcounter{equation}{0}

In this section we first present those versions of separation results for convex sets in vector spaces, which are expressed via their algebraic cores and are needed in what follows. This definitely has an overlapping with known separation theorems in vector spaces (see, e.g., \cite{hol}), while some statements and proofs are different. We show here that an ``extreme version" of the proper separation of a point from a convex set which is the core of itself, implies the standard Hahn-Banach extension theorem in vector spaces. Furthermore, the usage of proper separation allows us to derive a vector space algebraic core counterpart of a fundamental result of finite-dimensional convex geometry.

Recall that two nonempty subsets $\Omega_1$ and $\Omega_2$ of a vector space $X$  are {\em separated} by a hyperplane if there exists a nonzero linear function $f\colon X\to\R$ such that
\begin{equation}\label{convex separation0}
\sup\big\{f(x)\;\big|\;x\in\Omega_1\big\}\le\inf\big\{f(x)\;\big|\;x\in\Omega_2\big\}.
\end{equation}
If we have in addition that
\begin{equation}\label{convex proper separation}
\inf\big\{f(x)\;\big|\;x\in\Omega_1\big\}<\sup\big\{f(x)\;\big|\;x\in\Omega_2\big\},
\end{equation}
i.e., there exist vectors $x_1\in\Omega_1$ and $x_2\in\Omega_2$ with $f(x_1)<f(x_2)$, then the sets $\Omega_1$ and $\Omega_2$ are {\em properly separated} by a hyperplane. For brevity, we drop mentioning ``by a hyperplane" in what follows if no confusion arises.

The following observation shows that the separation notions in \eqref{convex separation0} and \eqref{convex proper separation} are equivalent to each other for the case where two sets $\Omega_1$ and $\Omega_2$ such that one of them is convex and its core is nonempty, while the other one is a singleton that does not belong to the set. \vspace*{-0.1in}

\begin{Proposition}\label{SP} Let $\Omega$ be a convex subset of $X$ with $\mbox{\rm core}(\Omega)\ne\emp$, and let $x_0\notin\Omega$. Then $\Omega$ and $\{x_0\}$ are separated if and only if they are properly separated.
\end{Proposition}\vspace*{-0.1in}
{\bf Proof.} It obviously suffices to show that if $\Omega$ and $\{x_0\}$ are separated, then they are properly separated as well. Choose $f\in X^\prime\setminus\{0\}$ such that
\begin{equation*}
f(x)\le f(x_0)\;\mbox{\rm for all }\;x\in\Omega.
\end{equation*}
Arguing by contradiction, suppose that for any $w\in\Omega$ we have $f(w)\ge f(x_0)$. This tells us that $f(x)=f(x_0)$ for all $x\in\Omega$, i.e., $f(x)$ is constant on $\Omega$. Since $\mbox{\rm core}(\Omega)\ne\emp$, Proposition~\ref{lm1} implies that the function $f(x)\equiv 0$ in $\Omega$, which cannot be true due the assumed separation of $\Omega$ and $\{x_0\}$. $\h$

Next we formulate the fundamental Hahn-Banach extension theorem in vector spaces; see, e.g., \cite[Theorem~I.6.A]{hol} for its proof. Recall that a function $p\colon X\to\R$ is {\em sublinear} if it is positively homogeneous and subadditive, i.e., $p(x_1+x_2)\le p(x_1)+p(x_2)$ for all $x_1,x_2\in X$. \vspace*{-0.3in}

\begin{Theorem}{\bf(Hahn-Banach theorem).}\label{HBT} Let $p\colon X\to\R$ be a sublinear function on $X$. Take a subspace $Y$ of $X$ and a linear function $g\colon Y\to\R$ satisfying
\begin{equation*}
g(y)\le p(y)\;\mbox{ whenever }\;y\in Y.
\end{equation*}
Then there exists a linear function $f\colon X\to\R$ such that $f(y)=g(y)$ for all $y\in Y$ and $f(x)\le p(x)$ for all $x\in X$.
\end{Theorem}\vspace*{-0.1in}

Given an absorbing set $\Omega$, define the {\em Minkowski gauge function} associated with $\Omega$ by
\begin{equation}\label{gauge}
p_\Omega(x):=\inf\big\{\lambda>0\;\big|\;x\in\lambda\Omega\big\}.
\end{equation}
In the case where $\Omega$ is convex, $p_\Omega\colon X\to \R$ is clearly sublinear on $X$.

Now we derive from the Hahn-Banach theorem \eqref{HBT} the basic proper separation result used in this paper. By the {\em extreme case} we understand the one where $\Omega=\mbox{\rm core}(\Omega)$. \vspace*{-0.1in}

\begin{Theorem}{\bf(proper separation theorem).}\label{CSL} Let $\Omega$ be a convex subset in $X$ with $\mbox{\rm core}(\Omega)\ne\emp$, and let $x_0\notin\Omega$. Then there exists a hyperplane that separates $\Omega$ and $\{x_0\}$ properly. In the case where $\Omega=\mbox{\rm core}(\Omega)$, there is a nonzero linear function $f\colon X\to\R$ with
\begin{equation}\label{opens}
f(x)<f(x_0)\;\mbox{ for all }\;x\in\Omega.
\end{equation}
\end{Theorem}\vspace*{-0.1in}
{\bf Proof.} Assume first that $0\in\mbox{\rm core}(\Omega)$, and so $\Omega$ is an absorbing set. Define the subspace $Y:=\mbox{\rm span}\{x_0\}$ and the function $g\colon Y\to\R$ by $g(\alpha x_0):=\alpha$ as $\alpha\in\R$. We intend to show that $g$ is linear and satisfies the estimate $g(y)\le p_\Omega(y)$ for all $y\in Y$ via the Minkowski gauge of $\Omega$ defined in \eqref{gauge}. Indeed, suppose that $y=\alpha x_0$ for some $\alpha\in\R$. If $\alpha\le 0$, then $g(y)=\alpha\le 0\le p_\Omega(y)$. If $\alpha>0$, then we get
$$
g(y)=\alpha\le\alpha p_\Omega(x_0)=p_\Omega(\alpha x_0)=p_\Omega(y).
$$
Since $p_\Omega$ is sublinear on $X$, the above Hahn-Banach theorem allows us to find a linear function $f\colon X\to\R$ such that $f(y)=g(y)$ for all $y\in Y$ and $f(x)\le p_\Omega(x)$ for all $x\in X$. The function $f$ is nonzero due to $f(x_0)=1$. This clearly yields
\begin{equation}\label{me}
f(x)\le p_\Omega(x)\le 1=f(x_0)\;\mbox{ for all }\;x\in\Omega,
\end{equation}
which justifies the separation property \eqref{convex separation0}. Proposition~\ref{SP} tells us that in fact we have the proper separation in this case.

Let us next examine the case where $0\notin\mbox{\rm core}(\Omega)$. Fix $a\in\mbox{\rm core}(\Omega)$ and consider the set $\Theta:=\Omega-a$ for which $0\in\mbox{\rm core}(\Theta)$. Then $\Theta$ and $\{x_0-a\}$ are property separated by the above, and thus $\Omega$ and $\{x_0\}$ are properly separated as well. Note finally that in the case where $\Omega=\mbox{\rm core}(\Omega)$ inequality \eqref{me} becomes strict, and hence we verify \eqref{opens}. $\h$

Now we present a characterization of the separation and proper separation for a singleton from a convex set that strengthens, in particular, the result of Theorem~\ref{CSL}.\vspace*{-0.1in}

\begin{Theorem}{\bf(characterization of proper separation of singletons from convex sets).}\label{lm_Sepa_co_point} Let $\Omega$ be a convex subset of $X$ with $\co(\Omega)\ne\emp$, and let $x_0\in X$. Then the following assertions are equivalent:\\[1ex]
{\rm(a)} $\Omega$ and $\{x_0\}$ are separated.\\[1ex]
{\rm(b)}  $\Omega$ and $\{x_0\}$ are properly separated.\\[1ex]
{\rm(c)} $x_0\notin\mbox{\rm core}(\Omega)$.
\end{Theorem}\vspace*{-0.1in}
{\bf Proof.} Recalling Proposition~\ref{SP}, it suffices to prove that (b) and (c) are equivalent. Firstly, suppose that $x_0$ and $\Omega$ are properly separated. Let $f\colon X\to\R$ be a nonzero linear function satisfying the condition
\begin{equation*}
f(x)\le f(x_0)\;\mbox{\rm for all }\;x\in\Omega,
\end{equation*}
and let the point $\ox\in\Omega$ satisfy the strict inequality
\begin{equation*}
f(\ox)<f(x_0).
\end{equation*}
Arguing by contradiction, suppose that $x_0\in \mbox{\rm core}(\Omega)$. Then we can choose $t>0$ such that $x_0+t(x_0-\ox)\in\Omega$. It tells us that
\begin{equation*}
f\big(x_0+t(x_0-\ox)\big)\le f(x_0)\;\mbox{\rm for all }\;x\in\Omega
\end{equation*}
and readily implies that $f(x_0)\le f(\ox)$, a contradiction.

To verify the converse statement, deduce from Propositions~\ref{core-conv} and \ref{2core} that $\co(\Omega)$ is a nonempty convex subset of $X$ with $\co(\co(\Omega))=\co(\Omega)\ne\emp$ and $x_0\notin\co(\Omega)$. Theorem~\ref{CSL} ensures that $x_0$ and $\co(\Omega)$ are properly separated, i.e., there exists a nonzero linear function $f\colon X\to\R$ such that
\begin{equation*}
f(x)\le f(\bar{x})\;\text{for all}\;x\in\co(\Omega),
\end{equation*}
and also there exists $w\in\co(\Omega)\subset\Omega$ such that $f(w) <f(x_0)$. Fix any $u\in\Omega$ and get by Proposition~\ref{aip} that $tw+(1-t)u\in\co(\Omega)$ whenever $0<t\le 1$. Then we have
\begin{equation*}
tf(w)+(1-t)f(u)=f\big(tw+(1-t)u\big)\le f(x_0).
\end{equation*}
Passing to the limit as $t\dn 0$ tells us that $f(u)\leq f(x_0)$, which verifies the proper separation of the point $x_0$ from the set $\Omega$. $\h$

Our next goal is to derive the Hahn-Banach theorem (Theorem~\ref{HBT}) from the extreme version of the proper separation result from Theorem~\ref{CSL}. Note the proof given below is different from the known relationships between the separation Hahn-Banach theorem, where the latter analytic result is derived by applying the full-scaled separation theorem to the epigraph and graph of the functions $p$ and $g$ given in Theorem~\ref{HBT}; see, e.g., \cite[Theorem~I.6.A]{hol}. To proceed, we first present the following lemma on sublinear functions.\vspace*{-0.1in}

\begin{Lemma}\label{lm-O-convex} Let $p\colon X\to\R$ be a sublinear function, and let $\Omega:=\{x\in X \;|\;p(x)<1\}$. Then the set $\Omega$ is convex and absorbing. Furthermore, we get that $p_\Omega=p$ for the Minkowski gauge function \eqref{gauge}, and that $\Omega=\mbox{\rm core}(\Omega)$.
\end{Lemma}\vspace*{-0.1in}
{\bf Proof.} The convexity of the set $\Omega$ obviously follows from the convexity of the sublinear function $p$. Since $\mbox{\rm core}(\Omega)\subset\Omega$, we only need to verify the opposite inclusion. Fix any $x_0\in\Omega$ and let $v\in X$ be arbitrary. If $p(v)=0$, then for any $0<\lambda<1$ we have
\begin{align*}
p(x_0+\lambda x)\le p(x_0)+\lambda p(v)=p(x_0)<1.
\end{align*}
In the case where $p(v)\ne 0$, define $\delta:=(1-p(x_0))/p(v)$ and observe that
\begin{align*}
p(x_0+\lambda v)&\le p(x_0)+\lambda p(v)\\
&<p(x_0)+\frac{1-p(x_0)}{p(v)}p(v)\\
&=p(x_0)+1-p(x_0)=1
\end{align*}
if $0<\lambda<\delta$. Thus we get $x_0+\lambda v\in\Omega$ for all such $\lambda$. It follows that $x_0\in\mbox{\rm core}(\Omega)$, and so $\text{core}(\Omega)=\Omega$. Observing that $0\in\Omega=\text{core}(\Omega)$, we see that the set $\Omega$ is absorbing.

Further, let us show that $p_\Omega=p$. Fix any $x\in X$ and check first that $p_\Omega(x)\le p(x)$. Picking $\lambda>p(x)$, we have $p(x/\lambda)<1$ implying that $x/\lambda\in\Omega$ and $x\in\lambda\Omega$. The definition of the Minkowski function tells us that $p_\Omega(x)\le\lambda$, and so $p_\Omega(x)\le p(x)$.

Finally, take $\lambda>0$ satisfying $x\in\lambda\Omega$. Then $x=\lambda w$ for some $w$, and hence $p(w)<1$. It shows that $p(x)=p(\lambda w)=\lambda p(w)<\lambda$, and so $p(x)\le p_\Omega(x)$. Since  $x\in X$ was chosen arbitrarily, we arrive at $p=p_\Omega$ and thus complete the proof. $\h$

Now we ready to derive the above Hahn-Banach theorem from the extreme case of proper convex separation in general vector spaces.\vspace*{-0.1in}

\begin{Theorem} {\bf(Hahn-Banach theorem follows from the extreme case of proper convex separation).} Let the result of Theorem~{\rm\ref{CSL}} hold for any convex subset $\Omega$ of a vector space $X$ with ${\rm core}(\Omega)\ne\emp$. Then we have the full statement of Theorem~{\rm\ref{HBT}}.
\end{Theorem}\vspace*{-0.1in}
{\bf Proof.} Fix in the framework of Theorem~\ref{HBT} a subspace $Y\subset X$, a linear function $g\colon Y\to\R$, and a sublinear function $p\colon X\to\R$ such that $g(y)\le p(y)$ for all $y\in Y$.
If $g=0$, then the zero function $f=0$ satisfies the requirements of the Hahn-Banach theorem. Thus it suffices to consider the case where  $g$ is nonzero. Then we can find $y_0\in Y$ with $g(y_0)=1$. Define the sets $\Omega:=\{x\in X\;|\;p(x)<1\}$ and $\Lambda:=\Omega+\ker g$. It follows from Lemma~\ref{lm-O-convex} that $\Omega$ is convex with $\mbox{\rm core}(\Omega)=\Omega$, which yields the convexity of $\Lambda$. Furthermore, we deduce from Proposition~\ref{lm-Core-equa} that $\mbox{\rm core}(\Lambda)=\Lambda$.

Observe next that $y_0\not\in\Lambda$. Indeed, suppose on the contrary that $y_0\in\Lambda$ and then get that $y_0=\omega+z$, where $p(\omega)<1$ and $g(z)=0$. It tells us that
\begin{equation*}
g(y_0)=g(\omega+z)=g(\omega)\le p(\omega)<1=g(y_0),
\end{equation*}
which is a contradiction. Using now the extreme case of Theorem~\ref{CSL} gives us a linear function $h\colon X\to\R$ such that
\begin{equation} \label{separation from lambda}
h(x)<h(y_0)\;\mbox{\rm for all }\;x\in\Lambda.
\end{equation}
Since $0\in\Lambda$, we have $0=h(0)<h(y_0)$. Define further a new linear function $f\colon X\to\R$ by
\begin{equation*}
f(x):=\frac{1}{h(y_0)}h(x)\;\mbox{ for all }\;x\in X.
\end{equation*}
We claim that $f$ is an extension of $g$ from $Y$ to $X$, and that $f(x)\le p(x)$ for all $x\in X$ as is stated in the Hahn-Banach theorem. To proceed, observe that $f(y_0)=1$ and verify that the inclusion $z\in\ker g$ (i.e., $g(z)=0$) implies that $f(z)=0$. By the contrary, suppose that $f(z)\ne 0$ and hence get that $h(z) = h(y_0)f(z)\ne 0$. It yields
\begin{equation*}
h\Big(\frac{h(y_0)}{h(z)}z\Big)=h(y_0),
\end{equation*}
which contradicts \eqref{separation from lambda} since $\frac{h(y_0)}{h(z)}z\in\ker f\subset\Lambda$. Thus we arrive at $f(z)=0$.

It is easy to see that $Y=\ker g\oplus\{y_0\}$, which allows us to find for any $y\in Y$ some $z\in\ker g$ and $\lambda\in\R$ such that $y=z+\lambda y_0$. Since $f(z)=0$ and $f(y_0)=1$, we have
\begin{align*}
f(y)=f(z+\lambda y_0)=f(z)+\lambda f(y_0)=\lambda=g(y),
\end{align*}
which clearly implies that $f\big|_Y=g$, i.e., $f$ is an extension of $g$ to the whole space $X$

To verify finally that $f(x)\le p(x)$ on $X$, pick $x\in X$ and fix a number $\lambda\ge 0$ with $x\in\lambda\Omega$, which is possible by the construction of $\Omega$. Having $x=\lambda\omega$ for some vector $\omega\in\Omega$, we deduce from the definition of $f$ that
\begin{align*}
f(x)=f(\lambda\omega)=\lambda f(\omega)=\frac{\lambda}{h(y_0)}h(\omega).
\end{align*}
Since $\omega\in\Omega\subset\Lambda$, it follows from \eqref{separation from lambda} that $h(\omega)<h(y_0)$ with $\frac{\lambda}{h(y_0)}\ge 0$, and so
\begin{align*}
f(x)=\frac{\lambda}{h(y_0)}h(\omega)\le\frac{\lambda}{h(y_0)}h(y_0)=\lambda.
\end{align*}
To complete the proof of the Hahn-Banach theorem, we obtain from the Minkowski gauge definition \eqref{gauge} and Lemma~\ref{lm-O-convex} that $f(x)\le p_\Omega(x)=p(x)$. $\h$

The final result of this section gives us a vector space counterpart of one of the most fundamental results of convex finite-dimensional geometry concerning relative interiors of convex sets. The following theorem is formulated similarly to \cite[Theorem~6.1]{r} with replacing the relative interior and the (topological) closure therein by the algebraic core and algebraic closure, respectively. The proof given below is based on the proper convex separation while being significantly different from the one in \cite{r} (see also \cite[Theorem~1.72]{bmn} for a modification with more details), which strongly exploits the finite-dimensional topology.\vspace*{-0.1in}

\begin{Theorem}\label{coreaffin} Let $\Omega$ be a convex subset of a vector space $X$. If $a\in\mbox{\rm core}(\Omega)$ and $b\in\mbox{\rm lin}(\Omega)$, then we have the interval inclusion $[a,b)\subset\mbox{\rm core}(\Omega)$.
\end{Theorem}\vspace*{-0.1in}
{\bf Proof.} Fix $\lambda\in(0,1)$, define $x_\lambda:=\lambda a+(1-\lambda)b$, and verify that $x_\lambda\in\mbox{\rm core}(\Omega)$. Arguing by contradiction, suppose that $x_{\lambda}\notin\co(\Omega)$. Then $\{x_{\lambda}\}$ and $\Omega$ can be properly separated by Theorem~\ref{lm_Sepa_co_point}. It means that there exists a nonzero linear function $f\colon X\to\R$ such that
\begin{equation}\label{eq_cor_lin}
f(x)\le f(x_{\lambda})=\lambda f(a)+(1-\lambda)f(b)\;\text{ for all }\;x\in\Omega.
\end{equation}
Since $b\in\mbox{\rm lin}(\Omega)$,  definition \eqref{lin} of the algebraic closure ensures the existence of $w\in\Omega$ such that $[w,b)\subset\Omega$. Thus for all natural numbers $n\in\mathbb{N}$ we have
\begin{equation*}
x_n:=b+\frac{1}{n}(w-b)\in\Omega.
\end{equation*}
This yields by \eqref{eq_cor_lin} the equivalence
\begin{equation*}
f(x_n)\le f(x_{\lambda})\Longleftrightarrow\frac{1}{n}f(w)-\frac{1}{n}f(b)+\lambda f(b)\le\lambda f(a),\quad n\in\mathbb{N}.
\end{equation*}
Passing to the limit a $n\to\infty$ gives us the inequality
\begin{equation}\label{eq-co-lin-im}
f(b)\le f(a).
\end{equation}
Since $a\in\co(\Omega)$, for any $m\in\mathbb{N}$ sufficiently large we have
\begin{equation*}
x_m:=a+\frac{1}{m}(a-b)\in\Omega
\end{equation*}
and then deduce from \eqref{eq_cor_lin} that
\begin{equation*}
f(x_m)=f(a)+\frac{1}{m}f(a)-\frac{1}{m}f(b)\le\lambda f(a)+(1-\lambda)f(b).
\end{equation*}
The passage there to the limit as $m\to\infty$ brings us to the equivalence
\begin{equation}\label{eq-co-lin-inver}
(1-\lambda)f(a)\le(1-\lambda)f(b)\Longleftrightarrow f(a)\le f(b)
\end{equation}
by $\lambda\in(0,1)$. Combining \eqref{eq-co-lin-im} and \eqref{eq-co-lin-inver}, we conclude that
\begin{equation}\label{eq_cor_lin-equa}
f(a)=f(b).
\end{equation}
It follows from $a\in\co(\Omega)$ that for any $v\in X$ there exits $t>0$ such that $a+tv\in\Omega$. Using finally \eqref{eq_cor_lin} and \eqref{eq_cor_lin-equa} tells us that $f(v)=0$. This is a contradiction, which verifies that $x_{\lambda}\in\Omega$ and thus completes the proof of the theorem. $\h$\vspace*{-0.1in}

\section{Algebraic Cores of Convex Graphs}\label{graph}
\setcounter{equation}{0}

This short section presents an extension of yet another important finite-dimensional result to the general framework of vector spaces. It concerns Rockafellar's theorem on representing relative interiors  of convex graphs of set-valued mappings via those for domain and image sets; see \cite[Theorem~6.8]{r} for an equivalent formulation. In \cite{cmn} we generalized this result to quasi-relative interiors of mappings between locally convex topological vector spaces under an additional quasi-regularity assumption that is always fulfilled in finite dimensions. The goal of this section is to derive a counterpart of the latter result for convex-graph mappings between arbitrary vector spaces in terms of algebraic cores of the involved sets without imposing any regularity assumptions. Our proof is based on the core characterization of proper separation of a point from a convex set that is given in Theorem~\ref{lm_Sepa_co_point}.

First we observe the following useful lemma.\vspace*{-0.1in}

\begin{Lemma}\label{Pro-Linear-Cor} Let $\Omega$ be a convex subset of a vector space $X$, and let $A\colon X\to Y$ be a linear operator. If $A$ is surjective $($i.e., $AX=Y$), then we have the inclusion
\begin{equation}\label{coreequality}
A\big(\mbox{\rm core}(\Omega)\big)\subset\mbox{\rm core}\big(A(\Omega)\big),
\end{equation}
which holds as equality if it is assumed in addition that $\mbox{\rm core}(\Omega)\ne\emp$.
\end{Lemma}\vspace*{-0.1in}
{\bf Proof.} Fix any $x_0\in\mbox{\rm core}(\Omega)$ and show that $A(x_0)\in\mbox{\rm core}(A(\Omega))$. Indeed, for every $v\in Y$ we have the surjectivity of $A$ that $v=A(u)$ with some $u\in X$. Choose $\delta>0$ such that $x_0+tu\in\Omega$ if $|t|<\delta$. Thus it follows that
\begin{equation*}
A(x_0)+tv=A(x_0+t u)\in A(\Omega)\;\mbox{\rm whenever }\;|t|<\delta,
\end{equation*}
which yields $A(x_0)\in\mbox{\rm core}(A(\Omega))$ and hence verifies \eqref{coreequality}. To prove the opposite inclusion
$$
\mbox{\rm core}\big(A(\Omega)\big)\subset A\big(\mbox{\rm core}(\Omega)\big),
$$
consider first the case where $0\in\mbox{\rm core}(\Omega)$. Choose any $y\in\mbox{\rm core}(A(\Omega))$ and find $t>0$ such that $y+ty\in A(\Omega)$, which tells us that
$$
y\in\frac{1}{1+t}A(\Omega)=A\Big(\frac{1}{1+t}\Omega\Big).
$$
Since $0\in\mbox{\rm core}(\Omega)$, it follows from Proposition~\ref{aip} that $\frac{1}{1+t}\Omega\subset\mbox{\rm core}(\Omega)$, and so $y\in A(\mbox{\rm core}(\Omega))$. It justifies the equality in \eqref{coreequality} in the case under consideration.

In the general case where $\mbox{\rm core}(\Omega)\ne\emp$, take any $a\in\mbox{\rm core}(\Omega)$ and get that $0\in\mbox{\rm core}(\Omega-a)$. It shows by the above that
$$
A\big(\mbox{\rm core}(\Omega-a)\big)=\mbox{\rm core}\big(A(\Omega-a)\big),
$$
which therefore verifies the equality in \eqref{coreequality} in the general case. $\h$

Now are are ready to establish the aforementioned vector space counterpart of Rockafellar's finite-dimensional theorem on convex graphs. Given a set-valued mapping $F\colon X\tto Y$ between vector spaces, its {\em domain} and {\em graph} are defined, respectively, by
\begin{equation*}
\dom(F)=\big\{x\in X\;\big|\;F(x)\ne\emp\}\;\;\mbox{ and }\;\;\gph(F):=\big\{(x,y)\in X\times Y\;\big|\;y\in F(x)\big\}.
\end{equation*}\vspace*{-0.4in}

\begin{Theorem}{\bf(algebraic cores of convex graphs in vector spaces).}\label{Theo-Rock-core} Let $F\colon X\tto Y$ be a convex set-valued mapping between vector spaces, and let $\mbox{\rm core}(\gph F)\ne\emp$. Then we have the following representation for the core of the graph:
\begin{equation}\label{core-gr}
\mbox{\rm core}(\gph F)=\big\{(x,y)\;\big|\;x\in\mbox{\rm core}(\dom F),\;y\in\mbox{\rm core}\big(F(x)\big)\big\}.
\end{equation}
\end{Theorem}\vspace*{-0.1in}
{\bf Proof.} Define the mapping $\mathcal{P}\colon X\times Y\to X$ by $(x,y)\mapsto x$. Then we clearly have that
\begin{equation*}
\mathcal{P}\big(\mbox{\rm core}\big(\gph F)\big)=\mbox{\rm core}\big(\mathcal{P}(\gph F)\big)=\mbox{\rm core}\big(\dom(F)\big).
\end{equation*}
It implies that $x_0\in\mbox{\rm core}(\dom(F))$ for any $(x_0,y_0)\in\mbox{\rm core}(\gph(F))$. In addition, for any $v\in Y$ there exists $\delta>0$ such that
\begin{equation*}
(x_0,y_0)+\lambda(0,v)\in\gph(F)\;\mbox{\rm whenever }\;|\lambda|<\delta.
\end{equation*}
It tells us that $y_0+\lambda v\in F(x_0)$ whenever $|\lambda|<\delta$, and so $y_0\in\mbox{\rm core}(F(x_0))$. This readily verifies the inclusion ``$\subset$" in \eqref{core-gr}.

To prove the opposite inclusion, fix any $(x_0,y_0)$ with $x_0\in \mbox{\rm core}(\dom(F))$ and $y_0\in\mbox{\rm core}(F(x_0))$. Arguing by contradiction, suppose that $(x_0,y_0)\notin\mbox{\rm core}(\gph F)$. By the proper separation result from Theorem~\ref{lm_Sepa_co_point} on the product space $X\times Y$, we find nonzero linear functions $f\colon X\to\R$ and $g\colon Y\to\R$ such that
\begin{equation*}
f(x)+g(y)\le f(x_0)+g(y_0)\;\mbox{\rm whenever }\;(x,y)\in\gph(F),
\end{equation*}
and furthermore there exists a pair $(\bar{x},\bar{y})\in \gph(F)$ satisfying
\begin{equation*}
f(\bar x)+g(\bar y)<f(x_0)+g(y_0).
\end{equation*}
If $x_0=\ox$, then we get the condition
\begin{equation*}
g(y)\le g(y_0)\;\mbox{\rm whenever }\;y\in F(x_0),
\end{equation*}
and thus $\bar{y}\in F(x_0)$ satisfies the strict inequality
\begin{equation*}
g(\bar y)<g(y_0).
\end{equation*}
The latter implies that $y_0\notin\mbox{\rm core}(F(x_0))$, a contradiction. It remains to consider the case where $x_0\ne\ox$. Then we can choose $t\in(0,1)$ to be so small that
$$
\tilde{x}:=x_0+t(x_0-\ox)\in\dom(F),
$$
which yields $x_0=\lambda\tilde{x}+(1-\lambda)\ox$ for some $0<\lambda<1$. Choosing further $\tilde{y}\in F(\tilde{x})$ gives us
\begin{equation}\label{1st}
f(\tilde{x})+g(\tilde{y})\le f(x_0)+g(y_0)
\end{equation}
and a pair $(\bar{x},\bar{y})\in\gph(F)$ satisfying
\begin{equation}\label{2nd}
f(\bar x)+g(\bar y)<f(x_0)+g(y_0).
\end{equation}
Multiplying \eqref{1st} by $\lambda$, \eqref{2nd} by $1-\lambda$, and then adding them together lead us to
\begin{equation*}
g(y^\prime)<g(y_0)\;\mbox{ with }\;y^\prime:=\lambda\tilde{y}+(1-\lambda)\oy\in F(x_0).
\end{equation*}
It yields $y_0\notin\mbox{\rm core}(F(x_0)$, a contradiction that completes the proof of the theorem. $\h$\vspace*{-0.1in}

\section{Normal Cone Intersection Rule via Algebraic Cores}\label{normal}
\setcounter{equation}{0}

In this section we begin the development of calculus rules for convex generalized differentiation in arbitrary vector spaces. Implementing the geometric approach to generalized differentiation, we start with normals to convex sets. Given a nonempty convex subset $\Omega$ of a vector space $X$, the {\em normal cone} to $\Omega$ at $\ox\in\Omega$ is defined by
\begin{equation}\label{nc}
N(\ox;\Omega):=\big\{f\in X^\prime\;\big|\;f(x-\ox)\le 0\;\text{ for all }\;x\in\Omega\big\}
\end{equation}
with the convention that $N(\ox;\Omega):=\emp$ if $\ox\notin\Omega$.

The main result of the normal cone calculus is the representation of the normal cone to convex set intersections via the normal cones to each set in the intersection. To derive such an {\em intersection rule}, we present first the following three lemmas concerning algebraic cores of convex sets. The first one gives us a simple formula for algebraic cores of set differences.\vspace*{-0.1in}

\begin{Lemma}\label{Co_sub_co} Let $\Omega_1$ and $\Omega_2$ be convex subsets in a vector space $X$ such that $\co(\Omega_1)\ne\emp$ and $\co(\Omega_2)\ne\emp$. Then we have the representation
\begin{equation*}
\co(\Omega_1-\Omega_2)=\co(\Omega_1)-\co(\Omega_2).
\end{equation*}
\end{Lemma}\vspace*{-0.1in}
{\bf Proof.} Define the linear mapping $A\colon X\times X\to X$ by $A(x,y):=x-y$ for all $(x,y)\in X\times X$. Then $A$ is a surjection. Letting $\Omega:=\Omega_1\times\Omega_2$ and using the product formula \eqref{prod} give us $\co(\Omega)=\co(\Omega_1)\times\co(\Omega_2)\ne\emp$. Applying now Lemma~\ref{Pro-Linear-Cor}, we have
\begin{equation*}
\co(\Omega_1-\Omega_2)=\co\big(A(\Omega)\big)=A\big(\co(\Omega)\big)=\co(\Omega_1)-\co(\Omega_2),
\end{equation*}
which completes the proof of the lemma. $\h$

The next lemma justifies proper convex separation of two convex sets expressed in terms of their algebraic cores. It is a direct consequence of the main separation result in Theorem~\ref{lm_Sepa_co_point} and the observation in the preceding lemma.\vspace*{-0.1in}

\begin{Lemma}\label{Theo_Sec_cone_ru} Let $\Omega_1$ and $\Omega_2$ be convex subsets of $X$ such that $\co(\Omega_1)\ne\emp$ and $\co(\Omega_2)\ne\emp$. Then the sets $\Omega_1$ and $\Omega_2$ are properly separated if and only if
\begin{equation}\label{qri-sep}
\co(\Omega_1)\cap\co(\Omega_2)=\emp.
\end{equation}
\end{Lemma}\vspace*{-0.1in}
{\bf Proof.} Define $\Omega:=\Omega_1-\Omega_2$ and get from Lemma~\ref{Co_sub_co} that condition \eqref{qri-sep} reduces to
\begin{equation*}
0\notin\co(\Omega_1-\Omega_2)=\co(\Omega_1)-\co(\Omega_2).
\end{equation*}
Theorem~\ref{lm_Sepa_co_point} tells us that the sets $\Omega$ and $\{0\}$ are properly separated, which clearly yields the proper separation of the sets $\Omega_1$ and $\Omega_2$.

To verify the opposite implication, suppose that $\Omega_1$ and $\Omega_2$ are properly separated. Then the sets $\Omega=\Omega_1-\Omega_2$ and $\{0\}$ are properly separated as well. By Theorem~\ref{lm_Sepa_co_point} we have
\begin{equation*}
0\notin\co(\Omega)=\co(\Omega_1-\Omega_2)=\co(\Omega_1)-\co(\Omega_2).
\end{equation*}
Thus $\co(\Omega_1)\cap\co(\Omega_2)=\emp$, which completes the proof. $\h$\vspace*{0.05in}

The last lemma in this section is a consequence of Theorem~\ref{Theo-Rock-core} that allows us to calculate algebraic cores of epigraphs for a special class of extended-real-valued functions. Recall that the {\em epigraph} of a function $\varphi\colon X\to\oR:=(-\infty,\infty]$ is given by
\begin{equation*}
\epi(\varphi):=\big\{(x,\alpha)\in X\times\R\;\big|\;\alpha\ge\varphi(x)\big\}.
\end{equation*}
\vspace*{-0.4in}

\begin{Lemma}\label{affine reg} Let $\Omega$ be a convex subset of a vector space $X$, and let $\co(\Omega)\ne\emp$. Given $f\in X^\prime$ and $b\in\R$, define the extended-real-valued function
\begin{equation*}
\psi(x):=\left\{\begin{aligned}
&f(x)+b\;&\mbox{ if }&\;x\in\Omega,\\
&\infty\;&\mbox{ if }\;&x\notin\Omega.
\end{aligned}\right.
\end{equation*}
Then we have the core representation for its epigraph:
\begin{equation*}\label{eq-core-epigrap}
\co\big(\epi(\psi)\big)=\big\{(x,\lambda)\in X\times\R\;\big|\;x\in\co(\Omega),\;\lambda>\psi(x)\big\}.
\end{equation*}
\end{Lemma}
{\bf Proof.} Let us first check that $\co\big(\epi(\psi)\big)\ne\emp$. Indeed, by $\co\big(\dom(\psi)\big)=\co(\Omega)\ne\emp$ there exists $\bar{x}\in\co\big(\dom(\psi)\big)$, and hence $(\bar{x},\bar{\lambda})=(\bar{x},\psi(\bar{x})+1)\in\epi(\psi)$. Taking any $(x,\lambda)\in X\times\R$, we show now that there exists $\delta>0$ such that
\begin{equation*}
(\bar{x},\bar{\lambda})+t(x,\lambda)\in\epi(\psi)\;\mbox{ whenever }\;|t|<\delta.
\end{equation*}
To proceed, we get from $\bar{x}\in\co(\Omega)$ a number $\delta_1>0$ ensuring that $\bar{x}+t x\in\Omega=\dom(\psi)$ for all $t$ with $|t|<\delta_1$. If $\lambda=f(x)$, then
\begin{equation*}
\psi(\bar{x}+t x)=f(\bar{x}+t x)+b\le\bar{\lambda}+t\lambda,
\end{equation*}
and hence $(\bar{x},\bar{\lambda})+t(x,\lambda)=(\bar{x}+t x,\bar{\lambda}+t\lambda)\in\epi(\psi)$ for all such $t$. In the remaining case where $\lambda\ne f(x)$, denote
\begin{equation*}
\delta:=\min\left\{\delta_1,\frac{1}{|\lambda-f(x)|}\right\}
\end{equation*}
and observe that for all $t$ with $|t|<\delta$ we have the equivalences
\begin{equation*}
\begin{aligned}&&f(\bar{x})+b-\bar{\lambda}=-1&\le t\big(\lambda-f(x)\big)\\
&\Longleftrightarrow& f(\bar{x}+ t x)+b&\le\bar{\lambda}+t\lambda\\
&\Longleftrightarrow&\psi(\bar{x}+t x)&\le\bar{\lambda}+t\lambda.
\end{aligned}
\end{equation*}
It means that $(\bar{x}+t x,\bar{\lambda}+t\lambda)\in\epi(\psi)$ for all $t$ with $|t|<\delta$. Therefore, we arrive at $(\bar{x},\bar{\lambda})\in\co(\epi(\psi))$, and hence get $\co\big(\epi(\psi)\big)\ne\emp$.

Define further the set-valued mapping $F\colon X\tto\R$ by $F(x):=[\psi(x),\infty)$. We easily see that $\dom(F)=\dom(\psi)=\Omega$ and $\gph(F)=\epi(\psi)$, which tells us that
\begin{equation*}
\co\big(\gph(F)\big)=\co\big(\epi(\psi)\big)\ne\emp.
\end{equation*}
Applying finally Theorem~\ref{Theo-Rock-core} verifies the conclusion of this lemma. $\h$

Now we are ready to establish the aforementioned normal intersection rule for finitely many sets in vector spaces under the core qualification condition.\vspace*{-0.1in}

\begin{Theorem}{\bf(normal cone intersection rule in vector spaces).}\label{Theointersecrule} Let $\Omega_1,\ldots,\Omega_m$ as $m\ge 2$ be convex subsets of a vector space $X$ under the qualification  condition
\begin{eqnarray}\label{eq1intersecrule}
\bigcap_{i=1}^m\co(\Omega_i)\ne\emp.
\end{eqnarray}
Then we have the normal cone intersection rule
\begin{equation}\label{eqintersecrule}
N\Big(\bar{x};\bigcap_{i=1}^m\Omega_i\Big)=\disp\sum_{i=1}^m N(\bar{x};\Omega_i)\;\mbox{ for all }\;\bar{x}\in\bigcap_{i=1}^m\Omega_i.
\end{equation}
\end{Theorem}\vspace*{-0.1in}
{\bf Proof.} We verify the claimed intersection rule for the case where $m=2$, while observing that the general case of finitely many sets can be easily deduced by induction. In fact, it suffices to prove the inclusion ``$\subset$" in \eqref{eqintersecrule} for $m=2$ by taking into account that the opposite inclusion is trivial. To proceed, fix $\ox\in\Omega_1\cap\Omega_2$ and $f\in N(\ox;\Omega_1\cap\Omega_2)$. Then the normal cone definition \eqref{nc} reads as
\begin{equation*}
f(x-\ox)\le 0\;\mbox{ for all }\;x\in\Omega_1\cap\Omega_2.
\end{equation*}
Consider further the auxiliary convex sets in the product space $X\times\R$ given by
\begin{equation}\label{Theta}
\Theta_1:=\Omega_1\times[0,\infty)\;\mbox{ and }\;\Theta_2:=\big\{(x,\lambda)\in X\times\R\;\big|\;x\in\Omega_2,\;\lambda\le f(x-\ox)\big\}.
\end{equation}
We deduce from \eqref{prod} that $\co(\Theta_1)=\co(\Omega_1)\times(0,\infty)$ and get by Lemma~\ref{affine reg} that
\begin{equation*}
\co(\Theta_2)=\big\{(x,\lambda)\in X\times\R\;\big|\;x\in\co(\Omega_2),\;\lambda< f(x-\ox)\big\}.
\end{equation*}
It obviously implies that $\co(\Theta_1)\cap\co(\Theta_2)=\emp$. Then the proper separation results of Lemma~\ref{Theo_Sec_cone_ru} applied to the sets in \eqref{Theta} gives us a nonzero pair $(h,\gamma)\in X^\prime\times\mathbb R$ such that
\begin{equation}\label{eq2tintersecrule}
h(x)+\lambda_1\gamma\le h(y)+\lambda_2\gamma\;\mbox{ for all }\;(x,\lambda_1)\in\Theta_1,\;(y,\lambda_2)\in\Theta_2,
\end{equation}
and that there exist pairs $(\tilde{x},\tilde{\lambda}_1)\in\Theta_1$ and $(\tilde{y},\tilde{\lambda}_2)\in\Theta_2$ satisfying
\begin{equation*}
h(\tilde{x})+\tilde{\lambda}_1\gamma< h(\tilde{y})+\tilde{\lambda}_2\gamma.
\end{equation*}
Observe that $\gamma\le 0$, since otherwise we get a contradiction by using \eqref{eq2tintersecrule} with $(\ox,1)\in\Theta_1$ and $(\ox,0)\in\Theta_2$. Now we employ the core qualification condition \eqref{eq1intersecrule} to show that $\gg\ne 0$. Suppose on the contrary that $\gamma=0$ and then get
\begin{equation*}
h(x)\le h(y)\;\mbox{ for all }\;x\in\Omega_1,\;y\in\Omega_2,\;\mbox{ and }\; h(\tilde{x})<h(\tilde{y})\;\mbox{ with }\;\tilde{x}\in\Omega_1,\;\tilde{y}\in\Omega_2.
\end{equation*}
This means that the original sets $\Omega_1$ and $\Omega_2$ are properly separated, and thus it follows from Lemma~\ref{Theo_Sec_cone_ru} that $\co(\Omega_1)\cap\co(\Omega_2)=\emp$, a contradiction showing us that $\gamma<0$.

Denoting further $\mu:=-\gamma>0$, we immediately deduce from \eqref{eq2tintersecrule} that
\begin{equation*}
 h(x)\le h(\ox)\;\mbox{ for all }\;x\in\Omega_1,\;\mbox{ and thus }\; h\in N(\ox;\Omega_1)\;\mbox{ and }\;\disp\frac{ h}{\mu}\in N(\ox;\Omega_1).
\end{equation*}
It also follows from \eqref{eq2tintersecrule}, due to $(\ox,0)\in\Theta_1$ and $(y,\alpha)\in\Theta_2$ with $\alpha:=f(y-\ox)$, that
\begin{equation*}
h(\ox)\le h(y)+\gamma f(y-\ox)\;\mbox{ whenever }\;y\in\Omega_2.
\end{equation*}
Dividing the both sides above by $\gamma$ and taking into account the linearity of the separating functions $f$ and $h$ ensure the inequality
\begin{equation*}
\Big(\frac{h}{\gamma}+f\Big)(y-\ox)\le 0\;\mbox{ for all }\;y\in\Omega_2,
\end{equation*}
and hence $\frac{ h}{\gamma}+f=-\frac{h}{\mu}+f\in N(\ox;\Omega_2)$. Therefore we arrive at
\begin{equation*}
 f\in\frac{ h}{\mu}+N(\ox;\Omega_2)\subset N(\ox;\Omega_1)+N(\ox;\Omega_2),
\end{equation*}
which verifies the claim in \eqref{eqintersecrule} for $m=2$, and thus completes the proof of the theorem. $\h$

As we see below, the normal cone intersection rule of Theorem~\ref{Theointersecrule} is crucial to derive major calculus rules for coderivatives and subgradients established in the subsequent sections.\vspace*{-0.1in}

\section{Coderivative and Subdifferential Sum Rules}\label{cod-sum}
\setcounter{equation}{0}

The main result of this section provides a sum rule for coderivatives of set-valued mappings and then uses it to derive the subdifferential sum rule for convex extended-real-valued functions on vector spaces. Note that the coderivative concept has never been considered in standard convex analysis; it came from variational analysis where it plays a prominent role. We refer the reader to the books \cite{m-book1,m-book,rw} and the bibliographies therein, where the reader can find important applications also to convex set-valued mappings between finite-dimensional and Banach spaces under closedness and lower semicontinuity assumptions.

It seems also that the coderivative notion has not been considered before in the setting of vector spaces without topology. Given a set-valued mapping $F\colon X\tto Y$ between arbitrary vector spaces $X$ and $Y$ and following the Banach space pattern, the {\em coderivative} of $F$ at $(\ox,\oy)\in\gph(F)$ is a set-valued mapping $D^*F(\ox,\oy)\colon Y'\tto X^\prime$ defined by
\begin{equation}\label{cod}
D^*F(\ox,\oy)(g):=\big\{f\in X^\prime\;\big|\;(f,-g)\in N\big((\ox,\oy);\gph(F)\big)\big\},\quad g\in Y'.
\end{equation}
Recall further that the {\em sum} of two set-valued mappings $F_1,F_2\colon X\tto Y$ is given by
\begin{equation*}
(F_1+F_2)(x)=F_1(x)+F_2(x):=\big\{y_1+y_2\in Y\;\big|\;y_1\in F_1(x),\;y_2\in F_2(x)\big\},\quad x\in X.
\end{equation*}
It is easy to see the domain relationship $\dom(F_1+F_2)=\dom(F_1)\cap\dom(F_2)$, and also that the graph of $F_1+F_2$ is convex provided that both mappings $F_1,F_2$ have this property. Our aim is to represent the coderivative of the sum $F_1+F_2$ at a given point of the graph in terms of the coderivatives of $F_1$ and $F_2$ at the corresponding points. We are going to derive a precise sum rule formula by using the normal cone intersection rule from Theorem~\ref{Theointersecrule}. To proceed, for any pair $(\ox,\oy)\in\gph(F_1+F_2)$ define the set
\begin{equation}\label{S}
S(\ox,\oy):=\big\{(\oy_1,\oy_2)\in Y\times Y\;\big|\;\oy=\oy_1+\oy_2,\;\oy_i\in F_i(\ox)\;\mbox{ as }\;i=1,2\big\}
\end{equation}
used in the formulation of the next theorem. Note that, in contrast to general coderivative sum rules in variational analysis and its convex specifications presented, e.g., in the aforementioned books, we do not impose now any uniform boundedness or inner semicompactness assumptions on \eqref{S} and also require a less restrictive qualification condition in comparison with that used in the variational analysis framework. The finite-dimensional version of the coderivative sum rule given below was first obtained in \cite[Theorem~11.1]{bmncal} via relative interiors, while its interior counterpart in linear convex topological vector spaces was established in \cite[Theorem~8.1]{mnrt}. Our new constraint qualification is expressed via algebraic cores in general vector spaces.\vspace*{-0.1in}

\begin{Theorem}{\bf(coderivative sum rule in vector spaces).}\label{CSR1} Let $F_1,F_2\colon X\tto Y$ be set-valued mappings with convex graphs between vector spaces, and let the following graphical core qualification condition be satisfied:
\begin{equation}\label{QCC}
\exists(x,y_1,y_2)\in X\times Y\times Y\;\mbox{\rm with }\;(x,y_1)\in\co\big(\gph(F_1)\big)\;\mbox{\rm and }\;(x, y_2)\in\co\big(\gph(F_2)\big).
\end{equation}
Then we have the coderivative sum rule
\begin{equation}\label{csr}
D^*(F_1+F_2)(\ox,\oy)(g)=D^*F_1(\ox,\oy_1)(g)+D^*F_2(\ox,\oy_2)(g)
\end{equation}
valid for all $(\ox,\oy)\in\gph(F_1+F_2)$, for all $g\in Y'$, and for all $(\oy_1,\oy_2)\in S(\ox,\oy)$.
\end{Theorem}\vspace*{-0.1in}
{\bf Proof.} Fix any $f\in D^*(F_1+F_2)(\ox,\oy)(g)$ and get by \eqref{cod} that $(f,-g)\in N((\ox,\oy);\gph(F_1+F_2))$. For every $(\oy_1,\oy_2)\in S(\ox,\oy)$ consider the convex sets
\begin{align*}
\Omega_1:=\big\{(x,y_1,y_2)\in X\times Y\times Y\;\big|\;y_1\in F_1(x)\big\},\;\Omega_2:=\big\{(x,y_1,y_2)\in X\times Y\times Y\;\big|\;y_2\in F_2(x)\big\}.
\end{align*}
It clearly follows from the constructions of $\Omega_i$, $i=1,2$, and the core definition \eqref{core} that
\begin{equation*}
\co(\Omega_i)=\big\{(x,y_1,y_2)\in X\times Y\times Y\;\big|\;(x,y_i)\in\co(\gph(F_i))\big\}\;\mbox{ for }\;i=1,2.
\end{equation*}
Furthermore, it is easy to deduce from the normal cone definition \eqref{nc} that
\begin{equation}\label{cod1}
(f,-g,-g)\in N\big((\ox,\oy_1,\oy_2);\Omega_1\cap\Omega_2\big).
\end{equation}
The imposed qualification condition \eqref{QCC} ensures that $\co(\Omega_1)\cap\co(\Omega_2)\ne\emp$. Applying now the intersection rule from Theorem~\ref{Theointersecrule} to the intersection in \eqref{cod1} brings us to
\begin{equation*}
(f,-g,-g)\in N\big((\ox,\oy_1,\oy_2);\Omega_1\big)+N\big((\ox,\oy_1,\oy_2);\Omega_2\big),
\end{equation*}
and therefore we arrive at the representation
\begin{equation*}
(f,-g,-g)=(f_1,-g,0)+(f_2,0,-g)\;\mbox{ with }\;(f_i,-g)\in N\big((\ox,\oy_i);\gph(F_i)\big),\;i=1,2,
\end{equation*}
The latter is is equivalent by the coderivative definition \eqref{cod} to
\begin{equation*}
f=f_1+f_2\in D^*F_1(\ox,\oy_1)(g)+D^*F_2(\ox,\oy_2)(g).
\end{equation*}
It readily justifies the inclusion ``$\subset$" in \eqref{csr}. The opposite inclusion is obvious. $\h$\vspace*{0.05in}

Next we present a direct consequence of Theorem~\ref{CSR1} to deriving the subdifferential sum rule for extended-real-valued convex functions on vector spaces under a new core qualification condition. The following result reduces to the classical one \cite[Theorem~23.8]{r} via the relative interior qualification condition in finite dimensions, while it does not require the continuity of one of the functions as in the known results in locally convex topological vector  spaces; see, e.g., the book \cite{z} and its references.

Considering an extended-real-valued convex function $\varphi\colon X\to\oR$ with the domain $\dom(\varphi):=\{x\in X\;|\;\ph(x)<\infty\}$ on a vector space $X$ and a point $\ox\in\dom\varphi$, an element $f\in X^\prime$ is a {\em subgradient} of $\varphi$ at $\ox$ if we have the inequality
\begin{equation*}
\varphi(x)\ge\varphi(\ox)+f(x-\ox)\;\mbox{ for all }\;x\in X.
\end{equation*}
As usual, the collection of all the subgradients of $\varphi$ at $\ox$ is called the {\em subdifferential} of $\varphi$ at $\ox$ and is denoted by $\partial\varphi(\ox)$. Remind that $\varphi$ is {\em proper} if $\dom(\varphi)\ne\emp$.\vspace*{-0.1in}

\begin{Theorem}{\bf(subdifferential sum rule in vector spaces).}\label{sr} Let $\varphi_i\colon X\to\oR$, $i=1,2$, be proper convex functions on a vector space $X$. Assume that the following epigraphical core qualification condition is satisfied:
\begin{equation}\label{riq}
\co\big(\epi(\varphi_1)\big)\cap\co\big(\epi(\varphi_2)\big)\ne\emp.
\end{equation}
Then we have the subdifferential sum rule
\begin{equation}\label{ssr}
\partial(\varphi_1+\varphi_2)(\ox)=\partial\varphi_1(\ox)+\partial\varphi_2(\ox)\;\mbox{ for all }\;\ox\in\dom(\varphi_1)\cap\dom(\varphi_2).
\end{equation}
\end{Theorem}\vspace*{-0.1in}
{\bf Proof.} Define the set-valued mappings $F_1,F_2\colon X\tto\R$ with convex graphs by
\begin{equation*}
F_i(x):=\big[\varphi_i(x),\infty\big)\;\mbox{ for }\;i=1,2.
\end{equation*}
Then the imposed qualification conditions \eqref{riq} tells us that the sets $\gph(F_i)=\epi(\varphi_i)$, $i=1,2$, have nonempty cores. Fixing any $\ox\in\dom(\varphi_1)\cap\dom(\varphi_2)$ and letting $\oy:=\varphi_1(\ox)+\varphi_2(\ox)$, we easily deduce from the definitions that
\begin{equation}\label{sr1}
f\in D^*(F_1+F_2)(\ox,\oy)(1)\;\mbox{ for every }\;f\in\partial(\varphi_1+\varphi_2)(\ox).
\end{equation}
Apply finally to \eqref{sr1} the coderivative sum rule from Theorem~\ref{CSR1} with $\oy_i=\varphi_i(\ox)$ as $i=1,2$. Then we arrive in this way at the relationships
\begin{equation*}
f\in D^*F_1(\ox,\oy_1)(1)+D^*F_2(\ox,\oy_2)(1)=\partial\varphi_1(\ox)+\partial\varphi_2(\ox),
\end{equation*}
which verify the inclusion ``$\subset$" in \eqref{ssr}. The opposite inclusion is trivial. $\h$\vspace*{-0.1in}

\section{Coderivative and Subdifferential Chain Rules}\label{cod-chain}
\setcounter{equation}{0}

This section deals with compositions of convex set-valued mappings between vector spaces and provides a precise chain rule to calculate coderivatives of compositions via coderivatives of their components under an appropriate core qualification condition for graphs. As a consequence of this general result, we derive a subdifferential chain rule for compositions of extended-real-valued convex functions and linear operators in the vector space setting.

Given two set-valued mappings $F\colon X\tto Y$ and $G\colon Y\tto Z$ between vector spaces, define their {\em composition} $(G\circ F)\colon X\tto Z$ by
\begin{equation*}
(G\circ F)(x)=\bigcup_{y\in F(x)}G(y):=\big\{z\in G(y)\;\big|\;y\in F(x)\big\},\quad x\in X,
\end{equation*}
and observe that $G\circ F$ is convex (i.e., its graph is convex) provided that both $F$ and $G$ have this property. Fix $\oz\in(G\circ F)(\ox)$ and consider the set
\begin{equation*}
M(\ox,\oz):=F(\ox)\cap G^{-1}(\oz).
\end{equation*}

The next theorem extends the finite-dimensional result of \cite[Theorem~11.2]{bmncal} expressed via relative interior (and the previous weaker versions of the coderivative chain rule discussed therein) to the general case of vector spaces with using the corresponding graphical core qualification condition. We refer the reader to \cite{m-book1,mnrt} for other infinite-dimensional coderivative chain rules under interior-type and related topological assumptions.\vspace*{-0.1in}

\begin{Theorem}{\bf(coderivative chain rule).}\label{scr} Let $F\colon X\tto Y$ and $G\colon Y\tto Z$ be convex set-valued mappings between vector spaces, and let there exist a triple $(x,y,z)\in X\times Y\times Z$ satisfying the graphical core qualification condition
\begin{equation}\label{QC1a}
(x,y)\in\co(\gph(F))\;\mbox{\rm and }\;(y,z)\in\mbox{\rm \co}(\gph(G)).
\end{equation}
 Then for any $(\ox,\oz)\in\gph(G\circ F)$ and $h\in Z'$ we have the coderivative chain rule
\begin{equation}\label{chain}
D^*(G\circ F)(\ox,\oz)(h)=D^*F(\ox,\oy)\circ D^*G(\oy,\oz)(h)\;\mbox{ whenever }\;\oy\in M(\ox,\oz).
\end{equation}
\end{Theorem}
{\bf Proof.} Fix $f\in D^*(G\circ F)(\ox,\oz)(h)$ and $\oy\in M(\ox,\oz)$. Then it follows from the coderivative definition that $(f,-h)\in N((\ox,\oz);\gph(G\circ F))$, which amounts to saying by \eqref{nc} that
\begin{equation*}
f(x-\ox)-h(z-\oz)\le 0\;\mbox{ for all }\;(x,z)\in\gph(G\circ F).
\end{equation*}
Define now the two convex sets in the product space $X\times Y\times Z$ by
\begin{equation*}
\Omega_1:=\gph(F)\times Z\;\mbox{ and }\;\Omega_2:=X\times\gph(G)
\end{equation*}
and easily deduce from these construction and and normal cone definition \eqref{nc} that
\begin{equation*}
(f,0,-h)\in N\big((\ox,\oy,\oz);\Omega_1\cap\Omega_2\big).
\end{equation*}
Then the core qualification condition \eqref{QC1a} tells us that $\co(\Omega_1)\cap\co(\Omega_2)\ne\emp$, and thus we are able to apply the normal cone intersection rule from
Theorem~\ref{Theointersecrule}. It leads us to
\begin{equation*}
(f,0,-h)\in N\big((\ox,\oy,\oz);\Omega_1\cap\Omega_2\big)=N\big((\ox,\oy,\oz);\Omega_1\big)+N\big((\ox,\oy,\oz);\Omega_2\big).
\end{equation*}
The latter yields the existence of an element $g\in Y'$ satisfying the representation $(f,0,-h)=(f,-g,0)+(0,g,-h)$ and the inclusions
\begin{equation*}
(f,-g)\in N\big((\ox,\oy);\gph(F)\big)\;\mbox{ and }\;(g,-h)\in N\big((\oy,\oz);\gph(G)\big).
\end{equation*}
Employing again the coderivative definition \eqref{cod} gives us the relationships
\begin{equation*}
f\in D^*F(\ox,\oy)(g)\;\mbox{ and }\;g\in D^*G(\oy,\oz)(h),
\end{equation*}
which justify the inclusion ``$\subset$" in \eqref{chain}. The opposite inclusion is straightforward. $\h$

As a simple consequence of Theorem~\ref{scr}, we present a subdifferential sum rule of convex analysis under an appropriate core qualification condition in vector spaces. It is a significant departure from the classical result of \cite[Theorem~23.9]{r} obtained under the relative interior qualification condition in finite dimensions, as well as from the corresponding infinite-dimensional chain rules given in topological frameworks; see, e.g., \cite{mnrt,z}.\vspace*{-0.1in}

\begin{Theorem}{\bf(subdifferential chain rule).}\label{sub-chain} Let $A\colon X\to Y$ be a linear mapping between vector spaces, and let $\varphi\colon Y\to\oR$ be a convex function. Assume that the range of $A$ contains a point of $\mbox{\rm\co(dom}(\varphi))$, and that $\co(\epi(\varphi))\ne\emp$. Then denoting $\oy:=A(\ox)\in\dom(\varphi)$ with some $\ox\in X$, we have the following subdifferential chain rule:
\begin{equation}\label{chain-r1}
\partial(\varphi\circ A)(\ox)=A^*\big(\partial\varphi(\oy)\big):=\big\{A^*g\;\big|\;g\in\partial \varphi(\oy)\big\},
\end{equation}
where $A^*\colon Y'\to X^\prime$ is the adjoint of a linear operator $A$ defined by
\begin{equation*}
A^*g(x):=g(Ax)\;\text{for all}\;g\in Y'\;\text{ and }\;x\in X.
\end{equation*}
\end{Theorem}
{\bf Proof.} Apply Theorem~\ref{scr} with $F(x):=\{A(x)\}$ and $G(x):=[\varphi(x),\infty)$. Then we have that $\co(\gph(F))=\co(\gph(A))$, $\co(\dom(G))=\co(\dom(\varphi))$, and $\gph(G)=\epi(\varphi)$. The imposed assumptions guarantee that the qualification condition \eqref{QC1a} is satisfied. Applying now Theorem~\ref{scr}, we get the equalities
\begin{equation*}
\partial(\varphi\circ A)(\ox)=D^*(G\circ A)(1)=D^*A\big(D^*G(\ox,\oy)(1)\big)=A^*\big(\partial\varphi(\oy)\big),
\end{equation*}
which verify \eqref{chain-r1} and hence completes the proof of the subdifferential chain rule. $\h$\vspace*{-0.1in}

\section{Subgradients of Convex Optimal Value Functions}\label{marginal}
\setcounter{equation}{0}

The last section of the paper is devoted to subdifferential study of the class of {\em optimal value/marginal functions} that are of high importance in variational analysis, optimization, and their numerous applications; see, e.g., \cite{bi,m-book1,m-book,rw} and the references therein. Extended-real-valued functions of this type are defined by
\begin{equation}\label{optimal value}
\mu(x):=\inf\big\{\ph(x,y)\;\big|\;y\in F(x)\big\},
\end{equation}
where $\ph\colon X\times Y\to\oR$, and where $F\colon X\tto Y$ is a set-valued mapping between vector spaces. Functions of type \eqref{optimal value} are intrinsically nondifferentiable even in the setting of smooth cost functions $\ph$ and simple moving sets $F(x)$. Clearly, such functions describe the optimal cost values in problems of parametric optimization
\begin{equation*}
\mbox{\rm minimize }\;\ph(x,y)\;\mbox{\rm subject to }\;y\in F(x),
\end{equation*}
but in fact the spectrum of their theoretical and applications is much broader; see the references above for more details and discussions. In particular, subdifferential information on $\mu(x)$ is crucial to understand behavior of marginal functions with respect to parameters.

Our goal here is to consider the case where $\ph$ and $F$ are convex in \eqref{optimal value}, and thus $\mu(x)$ is convex as well. We conduct our study in the framework of general vector spaces.

The case of convex subdifferentiation of functions \eqref{optimal value} is significantly different from nonconvex settings, where only upper estimates of subdifferentials are obtained under various qualification conditions; see \cite{m-book1,m-book,rw} with more references. As is known by now, the convex setting for \eqref{optimal value} allows us to derive a precise subdifferential formula for \eqref{optimal value} via the subdifferential of $\ph$ and the coderivative of $F$; see \cite{bmn}. To the best of our knowledge, the strongest result on calculating the convex subdifferential of \eqref{optimal value} in finite-dimensional spaces is obtained in \cite[Theorem~9.1]{bmncal} under a relative interior qualification condition. Its extension to locally convex topological vector spaces as given in \cite[Theorem~8.2]{mnrt} requires the continuity of $\ph$ in \eqref{optimal value} and does not reduce to \cite{bmn,bmncal} in finite dimensions. The following theorem is free of the continuity assumptions while imposing instead a much milder qualification condition in terms of algebraic cores of $\dom(\ph)$ and $\gph(F)$. It gives us back \cite[Theorem~9.1]{bmncal} when both spaces $X$ and $Y$ are finite-dimensional.\vspace*{-0.1in}

\begin{Theorem}{\bf(subdifferentiation of convex optimal value functions).}\label{mr3} Let $\mu(\cdot)$ be the optimal value function {\rm(\ref{optimal value})} generated by a convex mapping $F\colon X\tto Y$ between vector spaces and a convex extended-real-valued function $\ph\colon X\times Y\to\oR$. Suppose that $\co(\gph(F))\ne\emp$ and $\co(\epi(\ph))\ne\emp$, and that $\mu(x)>-\infty$ for all $x\in X$. Given $\ox\in\dom(\mu)$, consider the argminimum set
\begin{equation*}
S(\ox):=\big\{\oy\in F(\ox)\;\big|\;\mu(\ox)=\ph(\ox,\oy)\big\},
\end{equation*}
which is assumed to be nonempty. Then for any $\oy\in S(\ox)$ we have the equality
\begin{equation}\label{vf2}
\partial\mu(\ox)=\bigcup_{(f,g)\in\partial\ph(\ox,\oy)}\big[f+D^*F(\ox,\oy)(g)\big]
\end{equation}
provided that the following qualification condition is satisfied:
\begin{equation}\label{qf}
\co\big(\dom(\ph)\big)\cap\co\big(\gph(F)\big)\ne\emp.
\end{equation}
\end{Theorem}\vspace*{-0.1in}
{\bf Proof.} It is sufficient to verify the inclusion ``$\subset$" in \eqref{vf2}, since the opposite one is straightforward. To proceed, take any $h\in\partial\mu(\ox)$ and $\oy\in S(\ox)$ and then consider the sum
\begin{equation}\label{marg}
\psi(x,y):=\ph(x,y)+\delta_{{\rm\small gph}(F)}(x,y)\;\mbox{ for all }\;(x,y)\in X\times Y,
\end{equation}
where $\dd_\Omega(x)$ denotes the indicator function of a set $\Omega$ that equals $0$ if $x\in\Omega$ and $\infty$ otherwise.

Since the domain of $\delta_{{\rm\small gph}(F)}$ is $\gph(F)$ and the epigraph of $\delta_{{\rm\small gph}(F)}$ is $\gph(F)\times[0,\infty)$, it follows from \eqref{qf} and the application of Theorem~\ref{sr} to the summation function \eqref{marg} that
\begin{equation*}
(h,0)\in\partial\psi(\ox,\oy)=\partial\ph(\ox,\oy)+N\big((\ox,\oy);\gph(F)\big).
\end{equation*}
Thus we get from the above that
\begin{equation*}
(h,0)=(f_1,g_1)+(f_2,g_2)\;\mbox{ with }\;(f_1,g_1)\in\partial\ph(\ox,\oy)\;\mbox{ and }\;(f_2,g_2)\in N\big((\ox,\oy);\gph(F)\big),
\end{equation*}
which yields $g_2=-g_1$. Hence we arrive at the inclusion $(f_2,-f_1)\in N((\ox,\oy);\gph(F))$ meaning by definition \eqref{cod} that $f_2\in D^*F(\ox,\oy)(g_1)$. It tells us that
\begin{equation*}
h=f_1+f_2\in f_1+D^*F(\ox,\oy)(g_1),
\end{equation*}
which justifies the claimed inclusion``$\subset$" in (\ref{vf2}). $\h$\vspace*{-0.1in}

We provided a direct geometric approach to study convex generalized differentiation in vector spaces using qualification conditions based on the algebraic core. A similar approach would allow us to obtain similar calculus results in locally convex topological vector spaces based on the interior instead of the algebraic core. Conversely, it is possible to study convex generalized differentiation in locally convex topological vector spaces based on the interior and then obtain calculus results in vector spaces by using the strongest locally convex topology on the underlying space.

\end{document}